\documentclass{amsart}
\usepackage[usenames]{pstcol}
\usepackage{pst-node,pst-grad} 
\author{Valery Alexeev and Michel Brion}
\address{Department of Mathematics\\
  University of Georgia\\
  Athens, GA 30602, USA} \email{valery@math.uga.edu}
\address{Institut Fourier, B. P. 74\\
  38402 Saint--Martin d'H\`eres Cedex, France}
\email{Michel.Brion@ujf-grenoble.fr}

\title{Toric degenerations of spherical varieties}
\date{March 23, 2004}

\renewcommand{\thesubsection}%
{{\bf\S\  \arabic{section}.\arabic{subsection}}}

\newcommand{\myfrac}[2]{\genfrac\{\}{0cm}{0}{#1}{#2}}

\theoremstyle{plain}
\newtheorem{theorem}{Theorem}[section]

\newcommand{\Lie}{\operatorname{Lie}}
\newcommand{\LambdaR}{\Lambda_{\bR}}
\newcommand{\GT}{\operatorname{GT}}
\newcommand{\Qlam}{Q_{\uw_0}(\lambda)}
\newcommand{\Qw}{Q_{\uw_0}}
\newcommand{\Qwstd}{Q_{\uw_0^{\rm std}}}
\newcommand{\uwn}{\underline{w}_0}
\newcommand{\oM}{\overline{\operatorname{M}}}

\newcommand{\SP}{\operatorname{SP}}

\renewcommand{\ss}{\operatorname{ss}}

\newcommand{\Conv}{\operatorname{Conv}}

\newcommand{\GL}{\operatorname{GL}}
\newcommand{\SL}{\operatorname{SL}}

\newcommand{\Hom}{\operatorname{Hom}}

\newcommand{\gr}{\operatorname{gr}}
\newcommand{\bA}{{\mathbb A}}
\newcommand{\bB}{{\mathbb B}}
\newcommand{\bC}{{\mathbb C}}

\newcommand{\bK}{{\mathbb K}}

\newcommand{\bN}{{\mathbb N}}
\newcommand{\bP}{{\mathbb P}}
\newcommand{\bQ}{{\mathbb Q}}
\newcommand{\bR}{{\mathbb R}}
\newcommand{\bT}{{\mathbb T}}

\newcommand{\bZ}{{\mathbb Z}}

\newcommand{\cB}{{\mathcal B}}
\newcommand{\cC}{{\mathcal C}}

\newcommand{\cE}{{\mathcal E}}

\newcommand{\cO}{{\mathcal O}}

\newcommand{\cR}{{\mathcal R}}

\newcommand{\cX}{{\mathcal X}}

\newcommand{\inv}{^{-1}}

\newcommand{\uw}{\underline{w}}

\DeclareMathSymbol{\curvearrowright}{\mathrel}{AMSb}{"79}
\DeclareMathSymbol\rightsquigarrow {\mathrel}{AMSa}{"20}
\DeclareMathSymbol\square {\mathord}{AMSa}{"03}
\DeclareMathSymbol{\ltimes}         {\mathbin}{AMSb}{"6E}
\DeclareMathSymbol{\nmid}           {\mathrel}{AMSb}{"2D}
\DeclareMathSymbol{\twoheadrightarrow}  {\mathrel}{AMSa}{"10}

\newcommand{\ratmap}{- \kern -3pt \to}

\newcommand{\onto}{\twoheadrightarrow}

\newcommand{\Cone}{\operatorname{Cone}}

\newcommand{\supp}{\operatorname{supp}}

\newcommand{\Proj}{\operatorname{Proj}}

\newcommand{\Spec}{\operatorname{Spec}}

\theoremstyle{plain}

\newtheorem{lemma}[theorem]{Lemma}
\newtheorem{corollary}[theorem]{Corollary}

\newtheorem{proposition}[theorem]{Proposition}

\newtheorem{conjecture}[theorem]{Conjecture}

\theoremstyle{definition}

\newtheorem{definition}[theorem]{Definition}

\newtheorem{example}[theorem]{Example}

\newtheorem{remark}[theorem]{Remark}   
\newtheorem{remarks}[theorem]{Remarks}

\newtheorem*{acknowledgement}{Acknowledgment} 
\newtheorem{acknowledgements-no}{Acknowledgments} 

\theoremstyle{remark}

\begin{document}
\bibliographystyle{amsalpha}

\begin{abstract}
  We prove that any affine, resp.~polarized projective, spherical
  variety admits a flat degeneration to an affine, resp.~polarized
  projective, toric variety. Motivated by Mirror Symmetry, we give
  conditions for the limit toric variety to be a Gorenstein Fano, and
  provide many examples. We also provide an explanation for the limits
  as boundary points of the moduli space of stable pairs whose
  existence is predicted by the Minimal Model Program.
\end{abstract}

\maketitle

\section*{Introduction}

Let $G$ be a connected reductive group. We prove that any
affine, resp.~polarized projective, spherical $G$-variety admits a
flat degeneration to an affine, resp.~polarized projective, toric
variety. This is obtained as a corollary of a more general result
which applies to arbitrary, not necessarily spherical $G$-schemes. 
We construct several degenerations, labeled by reduced decompositions
$\uw_0$ of the longest element $w_0$ in the Weyl group of~$G$.

A number of special cases of this result have been proved before.
The case of flag and Schubert varieties was considered by many
authors, beginning with the work of Gonciulea and Lakshmibai
\cite{Gonciulea-Lakshmibai} for grassmanians and varieties of complete
flags, see also
\cite{Sturmfels_GrobnerBases,GrossbergZabcic,Chirivi,KoganMiller}.  

Caldero's paper \cite{Caldero} constructs several degenerations of
the flag variety $G/B$ that depend on the choice of the reduced
decomposition $\uw_0$. His paper contains an important ingredient
which we use for the general case. The choice of $\uw_0$ previously 
appeared in the study of degenerations of the Bott-Samelson
resolutions of $G/B$ by Grossberg and Karshon
\cite{GrossbergKarshon}. As tempting as this connection is, 
the two degenerations do not seem to be directly related, see the
discussion at the end of Section 6.

Regarding the more general case of spherical varieties, Kaveh
\cite{Kaveh_SAGBI} proves the existence of a toric degeneration in the
case where $G$ is the symplectic group $\SP_{2n}$.

We were led to the subject by trying to understand mirrors of
Calabi-Yau hypersurfaces in spherical varieties, following the model
suggested by Givental and others
\cite{Givental_Stationary,Batyrev-et-al}. In the latter paper, the
four authors argue that mirrors of Calabi-Yau hypersurfaces in a
partial flag variety $X=\SL_n/P$ can be found among Calabi-Yau
hypersurfaces in $X_0^*$, the toric variety mirror-dual to the
Gonciulea-Lakshmibai toric degeneration $X_0$ of $X$.

Therefore, in the projective case we pay special attention to the
question of when the toric limits are Gorenstein Fano varieties, i.e.,
when they correspond to reflexive polytopes. For an arbitrary
polarized spherical variety $(X,L)$, with toric degeneration
$(X_0,L_0)$, we describe the corresponding polytope $Q=Q_{\uw_0}(X,L)$.
An important feature of the general case is that $L_0$ may only be a
$\bQ$-line bundle, even if $X$ is Gorenstein Fano and $L=\cO(-K)$ 
is the anticanonical line bundle ; but we still have $L_0=\cO(-K)$ in
this case, i.e., $X_0$ is $\bQ$-Gorenstein Fano.

In the case where $X=G/P$ is a flag variety and $L=\cO(-K)$, we give a 
criterion for $X_0$ to be Gorenstein, i.e., for the polytope $Q$ to 
be reflexive. In particular, if $X = \SL_n/P$ there are
\emph{many} different reflexive polytopes corresponding to toric
degenerations and not just the Gelfand-Tsetlin polytope, which
corresponds to the Gonciulea-Lakshmibai degeneration. 

We conjecture that for groups of type $A_n$ the limits of 
$(G/B, \cO(-K))$ are Gorenstein. However, this fails for general types
as we show in Example~\ref{ex:E6} for a simple group $G$ of type $E_6$.

\bigskip

The paper is organized as follows. Section 1 introduces the necessary
technical tools: the dual canonical basis and its string
parameterizations, string cones and polytopes, Caldero's filtration.

Section 2 is devoted to the affine case. We construct limits of
an affine $G$-scheme $X$, which are toric if $X$ is spherical. We also
prove that these limits have rational singularities and trivial
canonical class, if this holds for $X$.

In Section 3, we adapt these results to the setting of polarized
projective varieties. In particular, we degenerate any polarized
spherical variety $(X,L)$ to a $\bQ$-polarized toric variety
$(X_0,L_0)$, which is a $\bQ$-Gorenstein Fano if so is $X$.
Further, we show that the moment polytope of $(X_0,L_0)$
projects onto the moment polytope of $(X,L)$, with fibers being string 
polytopes. The latter are the moment polytopes of toric limits of flag
varieties. 

In Section 4 we describe these string and moment polytopes and give a
criterion for their integrality, in particular, a criterion for the
anticanonical limit of a flag variety to be a Gorenstein Fano
variety. We prove that for the varieties $G/P$ associated with
minuscule or cominuscule weights, the limit toric varieties are
Gorenstein Fano. In general, we show that the string polytopes have a
number of special vertices where the associated toric variety is
Gorenstein.

Section 5 presents examples and counterexamples of string polytopes,
including the Gelfand-Tsetlin polytopes for classical groups. All of
these are integral in type $A$, but not in other types. We also
consider properties of string polytopes that do not depend on the
choice of $\uw_0$.

In Section 6, we explain how the above toric limits may be understood
geometrically as boundary points of a ``universal'' compact moduli
space of stable pairs, the existence of which is predicted by the
Minimal Model Program.

\begin{acknowledgement}
  The first author was partially supported by NSF grant No. 0101280.
\end{acknowledgement}

\section{Notation and background}

We begin by introducing notation that will be used throughout this
paper. We consider algebraic schemes, varieties and groups over the 
field $\bC$ of complex numbers. Let $G$ be a connected reductive
group, $B$ a Borel subgroup with unipotent radical $U$, and $T$ a
maximal torus of $B$, so that $B=TU$. Let $\Phi=\Phi(G,T)$ be the root
system of $(G,T)$, with subset of positive roots
$\Phi^+=\Phi(B,T)$. We denote by $\alpha_1,\ldots,\alpha_r$ the
corresponding simple roots, where $r$ is the semisimple rank of $G$.

Let $W$ be the Weyl group of $(G,T)$ and let $s_1,\ldots,s_r\in W$ be
the simple reflections associated with
$\alpha_1,\ldots,\alpha_r$. These generators of $W$ define the length
function $\ell$ on this group; let $w_0$ be the unique element of
maximal length. Then $N:=\ell(w_0)$ is the number of positive roots,
i.e., the dimension of the flag variety $G/B$.

Denote by $\Lambda$ the character group of $T$, also called the weight
lattice of $G$. Let $\Lambda^+$ be the subset of dominant weights; it
is stable under the involution 
\begin{displaymath}
\lambda\mapsto -w_0\lambda=:\lambda^*
\end{displaymath}
of $\Lambda$. There is a partial order on $\Lambda$, defined by
$\mu\le\lambda$ if there exist nonnegative integers $n_1,\ldots,n_r$
such that $\lambda-\mu=n_1\alpha_1+\cdots+n_r\alpha_r$. We put 
$\Lambda_{\bR}:=\Lambda\otimes_{\bZ} \bR$, then the convex cone
generated by $\Lambda^+$ in $\Lambda_{\bR}$ is the positive Weyl
chamber $\Lambda^+_{\bR}$. For $\lambda,\mu \in\Lambda^+$, we have 
$\mu\le\lambda$ if and only if: $\mu$ lies in the convex hull of the
Weyl group orbit $W\lambda$, and $\lambda-\mu$ lies in the root
lattice. In particular, there are only finitely many $\mu\in\Lambda^+$
that are smaller than a fixed $\lambda\in\Lambda^+$.

Recall that $\Lambda^+$ parametrizes the isomorphism classes of simple
(rational) $G$-modules, as follows. Consider the algebra 
\begin{displaymath}
A:=\bC[G]^U
\end{displaymath}
of regular functions on $G$, invariant under right multiplication by
$U$. Then $G\times T$ acts on $A$, where $G$ acts via left
multiplication, and $T$ acts on the right, since it normalizes $U$. We
have an isomorphism of $G\times T$-modules 
\begin{displaymath}
A \cong \bigoplus_{\lambda\in\Lambda^+} V(\lambda^*),
\end{displaymath}
where $V(\lambda^*)$ is a simple $G$-module of highest weight
$\lambda^*=-w_0\lambda,$ and $T$ acts on each $V(\lambda^*)$ via the
character $\lambda$. Thus, $V(\lambda^*)$ is the weight space
$A_\lambda$.

The assignement $\lambda\mapsto V(\lambda^*)$ is the desired
parametrization. Note that each $V(\lambda^*)\subset A$ has a
canonical linear form, $B$-eigenvector of weight $\lambda$, obtained
by restricting any function on $G/U$ to $B/U\cong T$. In other words,
$V(\lambda) = V(\lambda^*)^*$ has a canonical highest weight vector,
denoted by $v_{\lambda}$. The stabilizer in $G$ of the line spanned by
$v_{\lambda}$ is a parabolic subgroup containing $B$, denoted by
$P_\lambda$. Then $\lambda$ extends to a character of $P_\lambda$,
that we still denote by $\lambda$; let $L_\lambda$ be the
corresponding $G$-linearized line bundle on $G/P_\lambda$. We have
\begin{displaymath}
H^0(G/P_\lambda,L_\lambda) = V(\lambda^*).
\end{displaymath} 
Note that $P_\lambda=B$ if and only if $\lambda$ is regular, i.e., in
the interior of $\Lambda^+_{\bR}$.

We next recall results from \cite{Caldero}. The vector space $A$ has a
remarkable basis $\bB= (b_{\lambda,\varphi})$, where each
$b_{\lambda,\varphi}$ is an eigenvector of $T\times T\subseteq G\times
T$, of weight $\lambda$ for the right $T$-action. Specifically, for
fixed $\lambda$, the vectors $b_{\lambda,\varphi}$ form the dual basis
of the basis of $V(\lambda)$ consisting of the nonzero $bv_{\lambda}$,
where $b$ lies in the specialization at $q=1$ of Kashiwara-Lusztig's
\emph{canonical basis} $\cB$. We say that $\bB$ is the \emph{dual
  canonical basis} of $A$ (its elements are denoted by
$b^*_{\lambda,\varphi}$ in \cite{Caldero}, but we will drop the * as
we will only use the dual canonical basis). The parameters $\varphi$
are defined via the \emph{string parametrization} of $\bB$ (or $\cB$)
which depends on the choice of a reduced decomposition
\begin{displaymath}
\uw_0= (s_{i_1}, s_{i_2},\ldots, s_{i_N})
\end{displaymath}
of the longest element of $W$, i.e., 
$w_0 = s_{i_1} s_{i_2}\cdots s_{i_N}$. It allows to define an
injective map 
\begin{displaymath}
\iota_{\uw_0}:\bB \to \Lambda^+\times \bN^N,
~b_{\lambda,\varphi}\mapsto (\lambda,\varphi)=(\lambda,t_1,\ldots,t_N).
\end{displaymath}
The image of $\iota_{\uw_0}$ is the intersection of a rational convex
polyhedral cone $\cC_{\uw_0}$ of $\Lambda_{\bR}\times \bR^N$,
with the lattice $\Lambda\times\bZ^N$. The projection of $\cC_{\uw_0}$
to $\bR^N$ is the \emph{string cone} $C_{\uw_0}$, a rational
polyhedral convex cone of $\bR^N$. Define a map 
\begin{displaymath}
\pi:\Lambda_{\bR}\times\bR^N \to \bR^N, \quad
(\lambda,t_1,\ldots,t_N) \mapsto 
-\lambda + t_1 \alpha_{i_1} + \cdots + t_N \alpha_{i_N}.
\end{displaymath}
Then the weight of $b_{\lambda,\varphi}$ as a left $T$-eigenvector
equals $\pi(\lambda,t_1,\ldots,t_N)$.

A crucial feature of the dual canonical basis is the following
multiplicative property: 
\begin{displaymath}
b_{\lambda,\varphi} b_{\mu,\psi} = b_{\lambda+\mu,\varphi+\psi} 
+ \sum_{\gamma} c_{\lambda,\varphi;\mu,\psi}^{\gamma} \,
b_{\lambda+\mu,\gamma},
\end{displaymath}
where $c_{\lambda,\varphi;\mu,\psi}^{\gamma}= 0$ unless
$\gamma < \varphi + \psi$. Here $\le$ denotes the lexicographic
ordering on $\bN^N$. In particular, if a product 
$b_{\lambda,\varphi} b_{\mu,\psi}$ is a nonzero scalar multiple of 
a basis vector $b_{\nu,\gamma}$, then $\nu=\lambda + \mu$ and 
$\gamma = \varphi + \psi$.

{}From this, Caldero deduces the existence of an increasing filtration
of the algebra $A$ by $T\times T$-submodules, such that the associated
graded algebra, $\gr A$, is isomorphic to the algebra of the monoid
\begin{displaymath}
\Gamma_{\uw_0}: = i_{\uw_0}(\bB) = 
\cC_{\uw_0}\cap (\Lambda\times\bZ^N).
\end{displaymath}
In geometric terms, the affine variety 
\begin{displaymath}
G//U:=\Spec(A) 
\end{displaymath}
degenerates to
\begin{displaymath} 
(G//U)_0:=\Spec \bC[\Gamma_{\uw_0}]. 
\end{displaymath}
The latter is an affine toric
variety for the action of the torus $T\times \bT$, where we put
\begin{displaymath}
\bT:=(\bC^*)^N.
\end{displaymath} 
Further, the degeneration is compatible with the actions of $T\times
T$ on $G//U$ (regarding $T\times T$ as a subgroup of $G\times T$), and
on $(G//U)_0$ via the homomorphism of tori
\begin{displaymath}
T\times T\to T\times \bT,~
(x,y)\mapsto (x^{-1}y, \alpha_{i_1}(x),\ldots,\alpha_{i_N}(x)).
\end{displaymath}
Indeed, the induced homomorphism of character groups maps any 
$(\lambda,t_1,\ldots,t_N)\in \Lambda\times \bZ^N$ to 
$(-\lambda + t_1 \alpha_{i_1} + \cdots + t_N \alpha_{i_N},\lambda)
\in \Lambda\times\Lambda$. 

Actually, these results are obtained in \cite{Caldero} for a
simply-connected, semisimple group $G$, but their extension to an
arbitrary connected reductive group $G$ is immediate. Indeed, we have
$G=(G^{\ss}\times C)/Z$ where $G^{ss}$ is a simply-connected,
semisimple group, $C$ is a torus, and $Z$ is a finite central subgroup
of $G^{\ss}\times C$. Thus, $U\subset G^{\ss}$, the Weyl groups of $G$
and $G^{\ss}$ are identified, and $\bC[G]^U$ is the subalgebra of
$Z$-invariants in $\bC[G^{\ss}]^U\otimes \bC[C]$. The latter space has
a basis consisting of the $b_{\lambda,\varphi}\otimes \chi$, where 
$b_{\lambda,\varphi}$ lies in the dual canonical basis for
$G^{\ss}$, and $\chi$ is a character of $C$. Thus, a basis of
$\bC[G]^U$ consists of those $b_{\lambda,\varphi}\otimes \chi$
such that $\lambda+\chi$ restricts to the trivial character of $Z$.

The cone $\cC_{\uw_0}$ is described in \cite[\S1]{Littelmann} and
\cite[3.10]{Berenstein_Zelevinsky}:

\begin{theorem}\label{thm:cone_eqns}
  The cone $\cC_{\uw_0} \subset \Lambda_{\bR} \times \bR^N$ is the
  intersection of the preimage $\Lambda_{\bR}\times C_{\uw_0}$ of the
  string cone $C_{\uw_0} \subset \bR^N$ with the $N$ half-spaces
  \begin{displaymath}
    t_{k} \le \langle\lambda,  \alpha_{i_{k}}^{\vee}\rangle 
    - \sum_{\ell=k+1}^{N} 
    \langle \alpha_{i_\ell}, \alpha_{i_k}^{\vee}\rangle
    t_\ell,  \quad k=1, \dots, N.
  \end{displaymath}
\end{theorem}

Therefore, for every dominant weight $\lambda$ the basis
$(b_{\lambda,\varphi})$ is in bijection with the set of integral
points of the rational convex polytope $Q_{\uw_0}(\lambda)$ obtained
by intersecting the string cone with the cone defined by the above
$\lambda$-\emph{inequalities}; here and later, by integral points we
mean those in the lattice $\bZ^N$. Note that we may define
$Q_{\uw_0}(\lambda)$, more generally, for any
$\lambda\in\Lambda^+_{\bR}$. 

\begin{definition}\label{defn:string_polytope}
  We will call the polytopes $Q_{\uw_0}(\lambda)$ the 
  \emph{string polytopes}. 
\end{definition}

Let $\lambda\in\Lambda^+$. Under the projection 
\begin{displaymath}
\pi_{\lambda}:\bR^N \to \Lambda_{\bR}, \quad (t_1,\ldots,t_N)
\mapsto - \lambda + t_1 \alpha_{i_1}+\cdots + t_N \alpha_{i_N},
\end{displaymath} 
the integral points of $\Qlam$ map onto the weights of the $T$-module 
$V(\lambda^*)$. These weights are precisely those points $\mu$ in the
convex hull of the orbit $W\lambda^*=-W\lambda$, such that
$\lambda^*-\mu$ is in the root lattice. Applying this construction to
all positive integral multiples $n\lambda$ and using the equality
$Q_{\uw_0}(n\lambda)=nQ_{\uw_0}(\lambda)$, we see that in fact 
\begin{displaymath}
\pi_\lambda(\Qlam) =\Conv(W\lambda^*).
\end{displaymath}
Moreover, for each $\mu\in W\lambda^*$, the $\mu$-weight subspace of
$V(\lambda^*)$ is $1$-dimensional, so the point $\mu$ has exactly one
preimage in $\Qlam$. This preimage must be a vertex, since $\mu$ is a
vertex of $\Conv(W\lambda^*)$.

\begin{definition}\label{vertex}
  The preimages $\pi_\lambda\inv(\mu)$, $\mu\in W\lambda^*$, will
  be called the \emph{extremal weight vertices} of $\Qlam$ and denoted
  by $q_\mu$. 
\end{definition}

Note that each $q_\mu$ is integral. For example, the 
\emph{lowest weight vertex} $q_{-\lambda}$ is just the origin. If
$\lambda$ is regular, then the $\lambda$-inequalities are strict
inequalities at the origin, so that the tangent cone to 
the string polytope at $q_{-\lambda}$ is the whole string cone. 

On the other hand, the \emph{highest weight vertex} $q_{\lambda^*}$ is
the point where all the $\lambda$-inequalities are equalities. One
easily checks that
\begin{displaymath}
q_{\lambda^*} = (\langle \lambda^*,\beta_1^{\vee}\rangle,\ldots,
\langle \lambda^*,\beta_N^{\vee}\rangle),
\end{displaymath}
where 
 \begin{displaymath}
\beta_1 := \alpha_{i_1},\, \beta_2 : = s_{i_1}\alpha_{i_2},\,\ldots,\,
\beta_N := s_{i_1}\cdots s_{i_{N-1}}\alpha_{i_N}
\end{displaymath}
is the enumeration of the simple roots associated with the reduced
decomposition $\uw_0$.

Explicit inequalities for the string cone $C_{\uwn}$ (and therefore for
the string polytopes) in the case of groups of type $A_n$ are given by
Berenstein-Zelevinsky in \cite[Thm.3.14]{Berenstein_Zelevinsky}.
Littelmann \cite[Thm.4.2]{Littelmann} describes the cone $C_{\uwn}$
for the so called \emph{nice decompositions} $\uwn$, and this is
generalized in \cite[Thm.3.12]{Berenstein_Zelevinsky}. Also,
\cite[\S9]{Littelmann} contains formulas for a few particular
decompositions in types $F_4$ and $E_8$. 

By using \cite[Thm.4.2]{Littelmann}, one checks that the highest weight
vertex is in the interior of the string cone, whenever $\lambda$ is
regular and $\uw_0$ is a nice decomposition. Equivalently, the tangent
cone to the string polytope at $q_{\lambda^*}$ is defined by the
inequalities 
\begin{displaymath}
t_{k} + \sum_{\ell=k+1}^{N} 
    \langle \alpha_{i_\ell}, \alpha_{i_k}^{\vee}\rangle t_\ell \le 0,
\quad k=1,\ldots,N.
\end{displaymath}
In particular, this cone is generated by a basis of the lattice $\bZ^N$. 
We do not know if this extends to arbitary reduced decompositions and
regular dominant weights.

%%%%%%%%%%%%%%%%%%%%%%%%%%%%%%%%%%%%%%%%%%%%%%%%%%%%%%%%%%%%%%%%%%%%%%%%%%%
\section{Degenerations of affine $G$-algebras}
\label{sec:affine}

Let $R$ be an affine $G$-algebra, i.e., a finitely generated algebra
over $\bC$ on which $G$ acts rationally by algebra automorphisms. We
first recall results on the structure of the $G$-module $R$, for which
a general reference is \cite{Grosshans}. Let
$R^U$ be the subalgebra of $U$-invariants in $R$; this is an affine
$T$-algebra. Likewise, the algebra $R^G$ of $G$-invariants is finitely
generated, and each weight space $R^U_{\lambda}$ is a finitely
generated $R^G$-module; further, $R^U_0 = R^B = R^G$. We have
canonical isomorphisms of $R^G$-$G$-modules 
\begin{displaymath}
R \cong \bigoplus_{\lambda\in\Lambda^+}
\Hom^G(V(\lambda),R) \otimes V(\lambda) \cong 
\bigoplus_{\lambda\in\Lambda^+} R^U_{\lambda} \otimes V(\lambda),
\end{displaymath}
that map any 
$u\otimes v\in \Hom^G(V(\lambda),R) \otimes V(\lambda)$ to $u(v)\in R$, 
resp.~$u(v_{\lambda})\otimes v\in R^U_{\lambda} \otimes V(\lambda)$. 
The corresponding summands of $R$ are the \emph{isotypical components} 
$R_{(\lambda)}$. Their products in $R$ satisfy
\begin{displaymath}
R_{(\lambda)} R_{(\mu)} \subseteq 
\bigoplus_{\nu\in\Lambda^+,\ \nu\le\lambda+\mu}R_{(\nu)}.
\end{displaymath}

We now extend Caldero's filtration to this setting.
Choose a reduced decomposition $\uw_0$, then we obtain an isomorphism
of $R^G$-$T$-modules
\begin{displaymath} 
R \cong \bigoplus_{(\lambda,\varphi)\in\cC_{\uw_0}} 
R^U_{\lambda}\otimes b_{\lambda^*,\varphi},
\end{displaymath}
where $T$ acts on each space 
$R^U_{\lambda}\otimes b_{\lambda^*,t_1,\ldots,t_N}$ via the weight
$\pi_{\lambda}(t_1,\ldots,t_N)$.

Define a partial ordering $\le$ on $\Gamma_{\uw_0}$ by setting:
$(\mu,\psi) \le (\lambda,\varphi)$ if, either $\mu<\lambda$, or
$\mu=\lambda$ and $\psi\le \varphi$. Note that this ordering is
compatible with addition, and that there are only finitely many 
$(\mu,\psi)\in\Gamma_{\uw_0}$ that are smaller than a given
$(\lambda,\varphi)$. Now put
\begin{displaymath} 
R_{\le(\lambda,\varphi)} := 
\bigoplus_{(\mu,\psi)\le (\lambda,\varphi)} 
R^U_{\mu}\otimes b_{\mu^*,\psi}
\end{displaymath}
and define similarly $R_{<(\lambda,\varphi)}$. Both are finitely
generated $R^G$-submodules of $R$, stable under $T$. Further, the 
$R_{\le(\lambda,\varphi)}$ form an increasing, exhaustive filtration
of $R$, satisfying the multiplicative properties
\begin{displaymath}
R_{\le(\lambda,\varphi)} R_{\le(\mu,\psi)}\subseteq 
R_{\le(\lambda+\mu,\varphi+\psi)},
\end{displaymath}
\begin{displaymath}
R_{\le(\lambda,\varphi)} R_{<(\mu,\psi)}\subseteq 
R_{<(\lambda+\mu,\varphi+\psi)}.
\end{displaymath}
Indeed, we have 
\begin{displaymath}
(f_{\lambda}\otimes b_{\lambda^*,\varphi})
(f_{\mu}\otimes b_{\mu^*,\varphi})
\in (f_{\lambda}f_{\mu}\otimes b_{\lambda^*+\mu^*,\varphi+\psi})
+ R_{<(\lambda+\mu,\varphi+\psi)},
\end{displaymath}
as follows from the multiplicative properties of the isotypical
components and of the dual canonical basis. Thus, we may define the
associated graded algebra 
\begin{displaymath}
\gr R := \bigoplus_{(\lambda,\varphi)\in\Gamma_{\uw_0}}
R_{\le(\lambda,\varphi)}/R_{<(\lambda,\varphi)}.
\end{displaymath}
This is a $R^G$-$T$-algebra, graded by the monoid $\Gamma_{\uw_0}$; it
may be regarded as the space 
$\bigoplus_{(\lambda,\varphi)\in\Gamma_{\uw_0}} 
R^U_{\lambda}\otimes b_{\lambda^*,\varphi}$ 
with ``componentwise'' multiplication. This implies readily

\begin{proposition}\label{filtration}
With the preceding notation, the $R^G$-$T$-algebra $\gr R$ is
isomorphic to $(R^U \otimes \gr A)^T$, where $T$ acts on $R^U$
naturally, and on $\gr A$ via the inverse of its right action. As a
consequence, $\gr R$ is an affine $T\times\bT$-algebra, and its
$\bT$-invariant subring is isomorphic to $R^G$.
\end{proposition}

We would like to obtain $\gr R$ as the associated graded algebra of
$R$ for a $\bN$-filtration, instead of the multifiltration indexed by
the partially ordered monoid $\Gamma_{\uw_0}$. This is the content of
the following

\begin{proposition}\label{deformation}
There exist an affine $\bN$-graded $T$-algebra $\cR$ and a
$T$-invariant element $t\in\cR_1$ such that 

\begin{enumerate}
\item
$t$ is a nonzerodivisor in $\cR$, i.e., $\cR$ is flat over the
polynomial ring $\bC[t]$.

\item
The $\bC[t,t^{-1}]$-$T$-algebra $\cR[t^{-1}]$ is isomorphic to
$R[t,t^{-1}]$. 

\item 
The $T$-algebra $\cR/t\cR$ is isomorphic to $\gr R$.
\end{enumerate}
\end{proposition}

\begin{proof}
  We adapt the arguments of \cite[3.2]{Caldero} to our setting.
  By Proposition \ref{filtration}, the algebra $\gr R$ is finitely
  generated. Choose homogeneous generators
  $\bar{f}_1,\ldots,\bar{f}_n$ of respective degrees
  $(\lambda_1,\varphi_1)$, $\ldots$, $(\lambda_n,\varphi_n)$, and lift
  these generators to $f_1,\ldots, f_n\in R$. Then one easily checks
  that $f_1,\ldots,f_n$ generate the algebra $R$. Let $S$ denote the
  $\Gamma_{\uw_0}$-graded polynomial ring $\bC[x_1,\ldots,x_n]$, where
  $\deg(x_i):=(\lambda_i,\varphi_i)$. Choose homogeneous generators
  $\bar{g}_1,\ldots,\bar{g}_p$ of the kernel of the surjective
  homomorphism $S\to \gr R$, $x_i\mapsto \bar{f}_i$, and put
  $\deg(\bar{g}_i)=:(\mu_i,\psi_i)$. Then
  $\bar{g}_i(f_1,\ldots,f_n)\in R_{<(\mu_i,\psi_i)}$ for all $i$.
  Thus, we may find elements $g_1,\ldots,g_p$ of $S$ such that :
  $g_i(f_1,\ldots,f_n)=0$ and $g_i\in \bar{g}_i+S_{<(\mu_i,\psi_i)}$.
  It follows that the natural map $S/g_1S+\cdots +g_pS \to R$ is an
  isomorphism (since both sides are filtered algebras, and the
  associated graded map is an isomorphism).

Next consider all the differences $(\mu_i,\psi_i)-(\nu_j,\gamma_j)$,
where $(\nu_j,\gamma_j)$ is the degree of a homogeneous component of
$g_i-\bar{g}_i$. We claim that there exists a linear form
$e:\Lambda_{\bR}\times\bR^N\to \bR$ that takes positive integral
values at all nonzero points of $\Lambda^+\times\bN^N$, 
and at all these differences. Indeed, recall that either
$\mu_i>\nu_j$, or $\mu_i=\nu_j$ and $\psi_i>\gamma_j$. So we may take
$e= (A e_1,e_2)$, where $A$ is a large positive integer,
$e_1:\Lambda_{\bR}\to \bR$ is a linear form that takes positive
integral values at all positive roots, and $e_2:\bR^N\to \bR$ is a
linear form that takes positive integer values at all nonzero points
of $\bN^N$ and at all the (finitely many) points $\psi_i-\gamma_j$
such that $\lambda_i=\mu_j$ (the existence of $e_2$ follows from
\cite[Lemma 3.2]{Caldero} and its proof).

Now put for any nonnegative integer $m$:
\begin{displaymath}
R_{\le m}:=\langle f_1^{a_1}\cdots f_n^{a_n}~\vert~
e(\sum_{i=1}^n a_i(\mu_i,\psi_i))\le m\rangle\subset R.
\end{displaymath}
This defines a $T$-stable $\bN$-filtration of the algebra $R$, such
that the associated graded is our original $\gr R$. Let $t$ be an
indeterminate and consider the Rees algebra 
\begin{displaymath}
\cR:=\bigoplus_{m=0}^{\infty} R_{\le m} \, t^m \subset R[t].
\end{displaymath}
Then $\cR$ and $t$ have the desired properties; the finite generation
of $\cR$ follows from that of $\gr R$.
\end{proof}

(One easily checks that the Rees algebra associated with our
multifiltration of $R$ is not finitely generated. The reason is that
the set of nonnegative elements of $\Lambda\times\bZ^N$ spans a
convex cone which is not polyhedral; it is not even closed.)

Proposition \ref{deformation} may be formulated in geometric terms:
the affine $G$-algebra $R$ corresponds to an affine $G$-scheme
$X=\Spec(R)$, and $\gr R$, to an affine $T\times\bT$-scheme that
we denote by $X_0$. Now there exists a family of affine $T$-schemes 
\begin{displaymath}
\pi:\cX\to\bA^1 
\end{displaymath}
such that 

\begin{enumerate}
\item
$\pi$ is flat.

\item
$\pi$ is trivial with fiber $X$ over the complement of $0$ in
$\bA^1$.

\item
The fiber of $\pi$ at $0$ is isomorphic to $X_0$.
\end{enumerate}

We then say that $X$ \emph{degenerates to} $X_0$, and that $X_0$ 
\emph{is a limit of} $X$. Note that $X_0$ depends on the choice of a
reduced decomposition of $w_0$.

Next we consider an affine \emph{spherical} $G$-variety $X=\Spec(R)$,
i.e., $X$ is normal and contains a dense $B$-orbit. Equivalently, $R$
is a domain, and we have an isomorphism of $G$-modules
\begin{displaymath}
R\cong \bigoplus_{\lambda\in \Lambda^+\cap \Cone(X)} V(\lambda),
\end{displaymath}
where $\Cone(X)$ is a rational polyhedral convex cone in
$\Lambda^+_{\bR}$. In particular, the multiplicity of each
$V(\lambda)$ in $R$ is $0$ or $1$, so that $R$ is
\emph{multiplicity-free}. 
Then $\Cone(X)$ is uniquely determined by $X$; it is
called the \emph{weight cone} of $X$. Now Propositions \ref{filtration} 
and \ref{deformation} imply readily 

\begin{proposition}\label{spherical}
Let $X$ be an affine spherical $G$-variety with weight cone
$\Cone(X)$. Then $X$ degenerates to the affine toric 
$T\times \bT$-variety $X_0$ such that
\begin{displaymath}
\Cone(X_0)=(\Cone(X)\times\bR^N)\cap\cC_{\uw_0}.
\end{displaymath}
This degeneration is $T$-equivariant, where $T$ acts on $X_0$ via the
homomorphism $T\to T\times\bT$, 
$x\mapsto (x^{-1},\alpha_{i_1}(x),\ldots,\alpha_{i_N}(x))$.
\end{proposition}

Returning to an arbitrary $G$-scheme $X$ with limit $X_0$, many
geometric properties hold for $X$ if and only if they hold for $X_0$,
as can be shown along the lines of \cite[\S 6]{Popov}
(see also \cite[\S 18]{Grosshans}). We record three
such properties that will be of use in the sequel. 

For a normal variety $X$, we denote by $\cO(K_X)$ its canonical sheaf,
where $K_X$ is a canonical divisor; then the rational equivalence
class of $K_X$ is uniquely defined and called the 
\emph{canonical class} of $X$. Also, recall that $X$ 
\emph{has rational singularities}, if $R^if_*(\cO_Y) = 0$ for some
desingularization $f:Y\to X$ and for any $i\ge 1$; then $X$ is
Cohen-Macaulay with dualizing sheaf $\cO(K_X)$.

\begin{proposition}\label{singularities}
Let $X$ be an affine $G$-variety, with degeneration $X_0$
corresponding to the choice of a reduced decomposition of $w_0$. 
Then $X$ is normal (resp.~has rational singularities) if and only if
$X_0$ is normal (resp.~has rational singularities). Further, if
$mK_X \sim 0$ for some integer $m$, then $mK_{X_0}\sim 0$.
\end{proposition}

\begin{proof}
If $X_0$ is normal, then so is $X$ e.g. by \cite[Thm.18.3]{Grosshans}.
Conversely, if $X=\Spec(R)$ is normal, then so is $\Spec(R^U)$
by \cite[Thm.18.4]{Grosshans}. It follows easily that
$X_0=\Spec(R^U\otimes A)^T$ is normal as well.

The assertion on rational singularities is proved similarly, by using
the stability of rational singularities under deformations
\cite{Elkik} and under quotients by reductive groups \cite{Boutot}.

Now assume that $X$ (or $X_0$) is normal ; then $\cX$ is normal by
Serre's criterion. If $mK_X \sim 0$, then the sheaf $\cO(mK_{\cX})$
is isomorphic to $\cO_{\cX}$ outside $X_0$. Since both sheaves are
divisorial, and $X_0$ is irreducible, it follows that
$\cO(mK_{\cX})\cong \cO_{\cX}(nX_0)$ for some integer $n$. But $X_0$
is the divisor of the regular function $\pi$ on $\cX$, so that
$\cO(mK_{\cX})\cong\cO_{\cX}$. Restricting to $X_0$, we obtain
$\cO(mK_{X_0})\cong\cO_{X_0}$. Indeed, it suffices to check this on
the regular locus of $X_0$; the latter is contained in the regular
locus of $\cX$, where the adjunction formula 
$\cO(mK_{\cX})\vert_{X_0}\cong \cO(mK_{X_0})$ holds.
\end{proof}

\begin{remarks}
(i) Since toric varieties have rational singularities, Propositions
\ref{spherical} and \ref{singularities} yield another proof of
the rationality of singularities of spherical varieties
(\cite[Thm.10]{Popov}).

\noindent
(ii) For any reduced decomposition of $w_0$, the limit $(G//U)_0$ 
has trivial canonical class. Indeed, $G//U$ contains the homogeneous
space $G/U$ as an open orbit, with complement of codimension at least
$2$. Further, $K_{G/U}\sim 0$, since $U$ acts by
unimodular transformations on the tangent space to $G/U$ at the base
point. Thus, Proposition \ref{singularities} applies.

On the other hand, $(G//U)_0$ is singular if (say) the group $G$ is
simple of rank $\ge 2$. Indeed, $G//U$ is singular in that case.
\end{remarks}

%%%%%%%%%%%%%%%%%%%%%%%%%%%%%%%%%%%%%%%%%%%%%%%%%%%%%%%%%%%%%%%%%%%%%%

\section{Degenerations of polarized projective $G$-varieties}
\label{sec:projective}

\begin{definition}
A \emph{polarized (projective) $G$-variety} is a pair $(X,L)$, where
$X$ is a normal projective variety equipped with a $G$-action, and $L$
is an ample $G$-linearized invertible sheaf on $X$.
\end{definition}

For any integer $n$, the sheaf $L^n:=L^{\otimes n}$ is also
$G$-linearized. Thus, the space $H^0(X,L^n)$ is a 
finite-dimensional rational $G$-module. Further, the graded algebra
\begin{displaymath}
R(X,L):=\bigoplus_{n=0}^{\infty} H^0(X,L^n)
\end{displaymath}
is a finitely generated, integrally closed domain. We have:
$X=\Proj\, R(X,L)$, and $L^n =\cO_X(n)$ for all integers $n$.

Next we choose a reduced decomposition $\uw_0$ of $w_0$, which yields
the associated graded algebra, $\gr R(X,L)$. By Proposition
\ref{singularities}, this algebra is still a finitely generated,
integrally closed domain, with an action of 
$\bC^*\times T\times\bT$ such that the $\bC^*$-action defines a
positive grading. Thus, 
\begin{displaymath}
X_0 : = \Proj\, \gr R(X,L)
\end{displaymath}
is a projective $T\times\bT$-variety equipped with 
$T\times\bT$-linearized sheaves 
\begin{displaymath}
L_0^{(n)}:=\cO_{X_0}(n)
\end{displaymath}
for all integers $n$. By \cite[Thm.3.5]{Demazure}, there exists a 
$\bQ$-Weil divisor $D$ on $X_0$ such that 
$L_0^{(n)}\cong\cO_{X_0}(nD)$ for all $n$. (Here $\cO_{X_0}(nD)$
denotes the sheaf of rational functions $f$ on $X_0$ such that 
the $\bQ$-Weil divisor $(f) + n D$ has non-negative coefficients).
Further, $D$ is $\bQ$-Cartier and ample, i.e., the sheaf $L_0^{(m)}$
is invertible and ample for any sufficiently divisible integer $m>0$.
In particular, every sheaf $L_0^{(n)}$ is divisorial, i.e., it is the
sheaf of sections of an integral Weil divisor (called the floor of
$nD$).

We then say that $(X,L)$ \emph{degenerates to the $\bQ$-polarized
variety} $(X_0,L_0)$, and that $(X_0,L_0)$ \emph{is a limit of}
$(X,L)$. (To be completely correct, we should say that $(X,L^n)$
degenerates to $(X_0,L_0^{(n)})$ for all positive integers $n$, as
$L_0$ does not determine uniquely all the $L_0^{(n)}$.)
Since $\gr R(X,L)\cong R(X,L)$ as a $\bC^*\times T$-module,
each space $H^0(X_0,L_0^{(n)})$ is isomorphic to $H^0(X,L^n)$ as a
$T$-module, but it acquires a compatible action of the torus $\bT$.

Now Propositions \ref{deformation} and \ref{spherical} may be
generalized as follows.

\begin{theorem}\label{family}
Let $(X,L)$ be a polarized $G$-variety and choose a reduced
decomposition $\uw_0$. Then there exists a family of $T$-varieties
$\pi:\cX\to\bA^1$, where $\cX$ is a normal variety, together with
divisorial sheaves $\cO_{\cX}(n)$ ($n\in\bZ$), such that 

\begin{enumerate}
\item
$\pi$ is projective and flat.

\item
$\pi$ is trivial with fiber $X$ over the complement of $0$ in
$\bA^1$, and $\cO_{\cX}(n)\vert_X\cong L^n$ for all $n$.

\item
The fiber of $\pi$ at $0$ is isomorphic to $X_0$, and
$\cO_{\cX}(n)\vert_{X_0}\cong L_0^{(n)}$ for all $n$.
\end{enumerate}

If, in addition, $X$ is spherical, then $X_0$ is a toric variety under 
$T\times\bT$.
\end{theorem}

\begin{proof}
Let $R:=R(X,L)$ with filtration $(R_{\le m})$ and let 
$\cR:=\bigoplus_{m=0}^{\infty} R_{\le m}\, t^m$, as in the proof of
Proposition \ref{deformation}. Then the Rees algebra $\cR$ carries two
compatible $\bN$-gradings: by $m$, and by the degree inherited from that
of $R$. For the second grading, the subspace $\cR_0$ is just the
polynomial ring $\bC[t]$. Put $\cX:=\Proj \cR$ with corresponding
twisting sheaves $\cO_{\cX}(n)$. Then $\cX$ is equipped with a
morphism to $\bA^1$ such that the assertions (i), (ii) and (iii) are
easily checked.  

Since $R$ and $\gr R(X,L)$ are normal, it follows that $\cR$ is
normal, by Serre's criterion. So $\cX$ is normal as well.
Consider the variety $\hat\cX:=\Spec(\cR)$ and the morphism $t:\hat
\cX \to \bA^1$. The latter admits a section, corresponding to the
projection $\cR\to\cR_0=\bC[t]$. We still denote by 
$\bA^1\subset\hat\cX$ the image of this section. Now we have a morphism
\begin{displaymath}
p:\hat \cX - \bA^1 \to \cX
\end{displaymath}
which is a geometric quotient by the $\bC^*$-action corresponding to
the second grading. Further, using the normality of $\hat\cX$, one
checks that each $\cO_{\cX}(n)$ is the eigenspace of degree $n$ in the
sheaf $p_*\cO_{\hat\cX - \bA^1}$. This sheaf is reflexive, as 
$\hat \cX - \bA^1$ is normal and $p$ is equidimensional. Thus, all the 
$\cO_{\cX}(n)$ are reflexive as well. Since they have generic rank
$1$, it follows that they are divisorial.
\end{proof}

Next recall that polarized toric varieties correspond to integral convex
polytopes, and that Fano toric varieties correspond to reflexive
polytopes, see e.g. \cite{Batyrev94}. The assignement of a polytope to
a polarized toric variety may be generalized to arbitrary
polarized projective $G$-varieties by the following well-kown result.

\begin{lemma}\label{moment}
Let $(X,L)$ be a polarized $G$-variety. Then the 
points $\frac{\lambda}{n}\in\Lambda_{\bQ}$ such that:
$\lambda\in\Lambda^+$, $n$ is a positive integer, and the isotypical
component $H^0(X,L^n)_{(\lambda)}$ is nonzero, are exactly the 
rational points of a rational convex polytope $P(X,L)$ of
$\Lambda_{\bR}$. Further, $P(X,L^m) = m P(X,L)$ for any positive
integer $m$.
\end{lemma}

Indeed, $H^0(X,L^n)_{(\lambda)}\ne 0$ if and only if 
$H^0(X,L^{(n)})^U_{\lambda}\ne 0$. Further, the algebra $R(X,L)^U$
is finitely generated and $\Lambda^+\times\bN$-graded; let
$(f_i)$ be homogeneous generators, and $(\lambda_i,n_i)$ their
degrees. Then one easily checks that $P(X,L)$ is the rational convex
hull of the points $\frac{\lambda_i}{n_i}$.

\begin{definition}
$P(X,L)$ is the \emph{moment polytope} of the polarized $G$-variety
$(X,L)$.
\end{definition}

By positive homogeneity of the moment polytope, this definition
extends to $\bQ$-polarized varieties, in particular, to any limit
$(X_0,L_0)$ of $(X,L)$. We denote by $Q_{\uw_0}(X,L)$ the moment
polytope of that limit. It is a rational convex polytope in
$\Lambda_{\bR}\times\bR^N$, related to the moment polytope of $(X,L)$
by the following

\begin{theorem}\label{fibered}
  The first projection $p:\Lambda_{\bR}\times\bR^N \to \Lambda_{\bR}$
  restricts to a surjective map
\begin{displaymath}
  p:Q_{\uw_0}(X,L) \onto P(X,L), 
\end{displaymath}
with fiber over any $\lambda\in\Lambda^+_{\bR}$ being the string
polytope $Q_{\uw_0}(\lambda^*)$.
In particular, for $\lambda\in \Lambda^+$, the limit of the flag
variety $G/P_{\lambda}$ is a toric variety under $\bT$, and its moment
polytope is the string polytope $Q_{\uw_0}(\lambda^*)$.
\end{theorem}

\begin{proof}
  By definition, the rational points of $Q_{\uw_0}(X,L)$ are the pairs
  $(\frac{\lambda}{n},\frac{\varphi}{n})$ such that:
  $H^0(X,L^n)_{(\lambda)}\ne 0$, and
  $(\lambda^*,\varphi)\in\cC_{\uw_0}$. This implies the first
  assertion for rational $\lambda$, and hence for all $\lambda$
  since all involved polytopes are rational.
  
  Consider the $G$-variety $X=G/P_\lambda$ and the $G$-linearized line
  bundle $L=L_\lambda$. Then $L$ is ample, and
  $R(X,L)=\bigoplus_{n=0}^{\infty} \, V(n\lambda^*)$.  Thus, $X$ is
  spherical with moment polytope being the point $\lambda^*$. This
  implies the second assertion.
\end{proof}

\begin{remark}
The moment polytope $P(X,L)$ need not be integral. Consider, for
example, the group $\GL_n$ acting on the projectivization of the space
of $n\times n$ matrices by left multiplication. Then one checks that
the moment polytope for the line bundle $\cO(1)$ has nonintegral
vertices as soon as $n\ge 3$. 
\end{remark}

Finally, we consider limits of polarized spherical varieties which are
$\bQ$-Fano in the sense of the following

\begin{definition} We say that a normal projective variety $X$ is a 
\emph{$\bQ$-Gorenstein Fano} (or simply \emph{$\bQ$-Fano}) if its
anticanonical class $-K = -K_X$ is $\bQ$-Cartier and ample. If 
$-K_X$ is Cartier and ample, we say that $X$ is 
\emph{(Gorenstein) Fano}.
\end{definition}

\begin{theorem}\label{QFano}
Let $(X,L)$ be a polarized spherical $G$-variety with limit
$(X_0,L_0)$. If $X$ is a $\bQ$-Fano and $L=\cO(-mK_X)$ for some
positive integer $m$, then $X_0$ is a $\bQ$-Fano and
$L_0^{(n)}=\cO(-nmK_{X_0})$ for all $n$. 

Further, $L_0^{(n)}$ is invertible if and only if the polytope
$nQ(X,L)$ is integral; in this case, $nP(X,L)$ is integral as well.

In particular, if $X$ is a Fano and $L=\cO(-K_X)$, then
$L_0=\cO(-K_{X_0})$. Thus, $X_0$ is a Fano if and only if $Q(X,L)$ is
integral; in this case, $Q(X,L)$ is reflexive. 
\end{theorem}

\begin{proof}
The isomorphism $L_0^{(n)}=\cO(-nmK_{X_0})$ is proved by the argument
of Proposition \ref{singularities}. It implies, of course, that $X_0$
is $\bQ$-Fano.

If $nQ(X,L)$ is integral, then its vertices yield global sections of
$\cO_{X_0}(n)$, eigenvectors of the torus $\bT$. These sections have
no common zeroes (since they do not vanish simultaneously at any
$\bT$-fixed point), so that $\cO_{X_0}(n)$ is invertible. 
Conversely, if $\cO_{X_0}(n)$ is invertible, then $-nmK$ is an ample
(integral) Cartier divisor on the toric variety $X_0$. Thus, it is
generated by its global sections, and the preceding argument shows
that all the vertices of $nQ(X,L)$ are integral. This proves the
second assertion.
\end{proof}

\begin{remarks}
(i) For an arbitrary polarized $G$-variety $(X,L)$ with limit
$(X_0,L_0)$, the following conditions are equivalent: 

\noindent
(a) $X$ has rational singularities and $H^i(X,L^n)=0$ for
all $i\ge 1$ and $n\ge 0$.

\noindent
(b) $X_0$ has rational singularities and $H^i(X_0,L_0^{(n)})=0$ for
all $i\ge 1$ and $n\ge 0$.

Indeed, one checks that (a) (resp.~(b)) is equivalent to the
rationality of singularities of $\Spec R(X,L)$ 
(resp.~$\Spec R(X_0,L_0)$), and one applies Proposition
\ref{singularities}. 

\smallskip

\noindent
(ii) By the Kawamata-Viehweg vanishing theorem, the preceding
assumption (a) holds in the case where $X$ is a $\bQ$-Fano with log 
terminal singularities, and $L=\cO(-mK_X)$. Together with the
argument of Proposition \ref{singularities}, it follows that the limit 
$X_0$ is a $\bQ$-Fano, and $L_0^{(n)} = \cO(-mnK_{X_0})$ for all
$n$. This applies, for example, to Fano varieties with their
anticanonical polarization. 

\smallskip

\noindent
(iii) For a polarized spherical variety $(X,L)$, one has $H^i(X,L^n)=0$
for all $i\ge 1$ and $n\ge 0$ (as follows e.g. from the rationality of
singularities of $\Spec R(X,L)$). Thus, the Hilbert function 
$n\mapsto h^0(X,L^n)$ is a polynomial for $n\ge 0$, and coincides with
the Ehrhart function 
\begin{displaymath}
n\mapsto \#(Q_{\uw_0}(X,L^n)\cap \bZ^N) = 
\#(nQ_{\uw_0}(X,L)\cap \bZ^N)
\end{displaymath}
of the polytope $Q_{\uw_0}(X,L)$. This generalizes a result of
Okounkov \cite{Okounkov}.

\smallskip

\noindent
(iv) Let $X$ be a Fano spherical variety and $L=\cO(-K_X)$. If the
limit $X_0$ is nonsingular, then it is isomorphic to $X$ as a
$T$-variety. Indeed, $X_0$ is Fano by Theorem \ref{QFano}. Since $X_0$
is toric, it follows easily that $H^1(X_0,T_{X_0})=0$. This implies,
in turn, our rigidity statement by the argument of
\cite[Thm.3.1]{AlexeevBrion_Moduli}. 

As a consequence, any anticanonical limit of a flag variety $X=G/P$
is singular, except if $X$ is a product of projective spaces.

\end{remarks}

%%%%%%%%%%%%%%%%%%%%%%%%%%%%%%%%%%%%%%%%
\section{String and moment polytopes}
\label{sec:Resulting polytopes}

Recall that $p:\Lambda_{\bR}\times \bR^N\to \Lambda_{\bR}$
denotes the projection onto the first summand; it maps the cone
$\cC_{\uw_0}$ onto the positive Weyl chamber $\Lambda^+_{\bR}$. 

\begin{definition}
  For every $\lambda\in \LambdaR$, let $\sigma^0_{\lambda}$ be the
  intersection of the images $p(\tau^0)$, where $\tau$ is a face of
  the cone $\cC_{\uw_0}$ with relative interior $\tau^0$, and
  $p(\tau^0)$ contains $\lambda$. Let $\Sigma_{\uw_0}$ be the
  collection of the closed cones
  $\sigma_\lambda:=\overline{\sigma^0_\lambda}$. 
\end{definition}

Since the intersection of any two faces of $\cC_{\uwn}$ is a face,
the intersection of any two cones of $\Sigma_{\uw_0}$ is again in
$\Sigma_{\uw_0}$. Hence, $\Sigma_{\uw_0}$ is a fan, with support
$\Lambda^+_{\bR}$. The following Lemma follows by an elementary
convexity argument.

\begin{lemma}
  Two weights $\lambda, \mu \in \Lambda^+_{\bR}$ lie in the
  same cone of $\Sigma_{\uw_0}$ if and only if the
  polytope $Q_{\uw_0}(\lambda+\mu)$ is the Minkowski sum of the 
  polytopes $Q_{\uw_0}(\lambda)$ and $Q_{\uw_0}(\mu)$.
\end{lemma}

\begin{corollary}\label{cor:string_fan}
  If the fan $\Sigma_{\uw_0}$ is trivial, i.e. consists of the
  positive chamber $\Lambda^+_{\bR}$ and its faces, then for all
  regular dominant weights $\lambda$ the polytopes
  $Q_{\uw_0}(\lambda)$ have the same normal fan. In this case, the 
  polytope $Q_{\uw_0}(X,L)$ is integral if and only if the string
  polytope $Q_{\uw_0}(\lambda)$ is integral for any vertex $\lambda$
  of $P(X,L)$.
\end{corollary}

If the fan $\Sigma_{\uw_0}$ is non-trivial then clearly the vertices
of $Q_{\uw_0}(X,L)$ are among the vertices of the polytopes
$Q_{\uw_0}(\lambda)$ as $\lambda$ goes over the vertices
in the subdivion $\Sigma_{\uw_0} \cap P(X,L)$ of $P(X,L)$.

By Theorem \ref{QFano}, the anticanonical limit of the flag variety
$G/B$ is a Fano if and only if the polytope 
$Q_{\uw_0}(G/B, \cO(-K)) = Q_{\uwn}(2\rho)$ is integral, where $\rho$ 
is the half-sum of positive roots. This is not true in general, as we
show below in Example~\ref{ex:E6}. However, we will show that any
toric limit $(X_0,L_0)$ of $(G/P_\lambda,L_\lambda)$ is Gorenstein at
its $\bT$-fixed points associated with extremal weight vertices (in
the sense of Definition \ref{vertex}). Recall that any such vertex
$q_\mu$ is integral, and hence defines a global section $s_\mu$ of
$L_0$, eigenvector of $\bT$ of weight $v_\mu$. The complement of the
zero locus of $s_\mu$ is a $\bT$-stable affine open subset of $X_0$,
denoted by $X_{0,\mu}$.

\begin{theorem}\label{extremal}
  For every extremal weight vertex $q_\mu$ of $\Qlam$, the
  corresponding open subset $X_{0,\mu}$ has trivial canonical class.
\end{theorem}

\begin{proof}
  Let $\pi:\cX\to\bA^1$ and $\cO_{\cX}(1)$ be as in Proposition
  \ref{family}. The extremal weight vertex $q_\mu=\pi_\lambda\inv(\mu)$
  defines a $1$-dimensional submodule of $\pi_*\cO_{\cX}(1)$ and hence
  a global section $s_\mu \in H^0(\cX,\cO_{\cX}(1))$. The family
  $\cX_\mu=\{s_\mu\ne 0\}\to\bA^1$ is an affine family whose fiber
  at any nonzero point is the complement of the zero locus in
  $G/P_\lambda$ of the extremal weight vector 
  $v_\mu \in H^0(G/P_\lambda,L_\lambda)$. This fiber is a
  $W$-translate of the negative cell $U^-P_\lambda/P_\lambda$; in
  particular, it is an affine space. Hence, it has trivial canonical
  class, so that the same holds for $X_{0,\mu}$ by the proof of
  Proposition \ref{singularities}.
\end{proof}

Next recall that a dominant weight $\lambda$ is called
\emph{minuscule} if the only weights of $V(\lambda)$ are the
extremal weights $\mu\in W\lambda$; equivalently, the only integral
points in $Q_{\uw_0}(\lambda)$ are the extremal weight vertices. 
Clearly, the set of minuscule weights is stable under
$\lambda\mapsto\lambda^*$. All the minuscule weights are fundamental; 
their list is as follows:

\begin{itemize}
\item Type $A_n$: any fundamental weight.
\item Type $B_n$: $\omega_n$.
\item Type $C_n$: $\omega_1$.
\item Type $D_n$: $\omega_1,\omega_{n-1},\omega_n$.
\item Type $E_6$: $\omega_1,\omega_6$.
\item Type $E_7$: $\omega_7$. 
\item Types $G_2$, $F_4$, $E_8$: none.
\end{itemize}

Also recall that a fundamental weight $\lambda$ is called
\emph{cominuscule} if the corresponding fundamental weight in the dual
root system is minuscule. The cominuscule weights are, of course, the
minuscule weights in type $A$, $D$ and $E$, plus the following:

\begin{itemize}
\item Type $B_n$: $\omega_1$.
\item Type $C_n$: $\omega_n$.
\end{itemize}

Thus, the flag varieties $G/P_\lambda$ associated with the
(co)minuscule weights are the grassmanians, the quadrics,
the grassmanians of maximal isotropic subspaces (with respect to a
symmetric or alternating non degenerate bilinear form), and two
exceptional varieties.

\begin{theorem}\label{minuscule}
  For any (co)minuscule weight $\lambda$ and for any reduced
  decomposition $\uw_0$, the polytope $Q_{\uw_0}(\lambda)$ is
  integral. Further, the corresponding toric variety $(G/P_\lambda)_0$
  is a Fano. In the minuscule case, all the vertices of
  $Q_{\uw_0}(\lambda)$ are extremal weight vertices.
\end{theorem}

\begin{proof}
We begin with the minuscule case. For any positive integer $n$, we
show that the dual canonical basis of $V(n\lambda^*)$ coincides (up to
nonzero scalar multiples) with its \emph{standard monomial basis}. The
latter can be defined as follows (see e.g. \cite{Gonciulea-Lakshmibai}). 
Let $W^{\lambda}$ be the set of representatives of minimal length of
all cosets $wW_\lambda$, where $W_\lambda$ is the isotropy group of
$\lambda$ in $W$. Then we may index the extremal weight vectors of
$V(\lambda^*)$ by $W^\lambda$, via 
\begin{displaymath}
x\mapsto p_x:=v_{-x(\lambda)}\in 
H^0(G/P_\lambda,L_{\lambda})=V(\lambda^*).
\end{displaymath}
The set $W^\lambda$ is partially ordered by the Bruhat ordering of
$W$. Now the products 
\begin{displaymath}
p_{x_1}\cdots p_{x_n} \in 
H^0(G/P_\lambda,L_{\lambda}^n)=V(n\lambda^*),
\end{displaymath}
where $x_1\ge \cdots\ge x_n$ in $W^\lambda$, form a basis of
$V(n\lambda^*)$. These products are called the standard monomials.
Further, if $x\ge y\in W^\lambda$, then there exists
a sequence $(s_{i_1},\ldots,s_{i_\ell})$ of simple reflections such
that 
\begin{displaymath}
x = s_{i_1}\cdots s_{i_\ell}y \quad \text{and} \quad 
\ell(x) = \ell + \ell(y),
\end{displaymath}
as follows from \cite[Lem.7.1]{Gonciulea-Lakshmibai}; we say that the
Bruhat ordering of $W^{\lambda}$ is generated by simple reflections. 
Thus, given $x_1\ge \cdots\ge x_n$ in $W^\lambda$, there exists a
reduced decomposition $w_0=s_{i_1}\cdots s_{i_N}$ such that each $x_k$
is a subword $s_{i_1}\cdots s_{i_j}$ for some $j=j(k)$. Now we may
apply \cite[Thm.14]{Caldero-Littelmann} which asserts that a nonzero
scalar multiple of $p_{x_1}\cdots p_{x_n}$ lies in the dual canonical
basis of $V(n\lambda^*)$.

Together with the multiplicative property of this basis, it follows
that any integral point of $Q_{\uw_0}(n\lambda)$ is a sum of $n$
extremal weight vertices. Since
$Q_{\uw_0}(n\lambda)=n\Qw(\lambda)$, it follows in turn that the
polytope $Q_{\uw_0}(\lambda)$ is integral and that its integral points
are precisely the extremal weight vertices.

Now, if $\lambda$ is cominuscule, then it is of classical type in the
sense of \cite{Lakshmibai-Seshadri}. Then a basis for $V(\lambda^*)$
consists of $T$-eigenvectors $p_{x,y}$, where $x$, $y\in W^{\lambda}$
form an admissible pair in the sense of [loc.cit.]; in particular, 
$x\ge y$ (and $x=y$ yields the extremal weight vectors). Further, by
[loc.cit.] again, a basis of $V(n\lambda^*)$ consists of the standard
monomials of degree $n$ in the $p_{x,y}$, i.e., of those products
$p_{x_1,y_1}\cdots p_{x_n,y_n}$ such that 
$x_1\ge y_1\ge\cdots\ge x_n\ge y_n$. Since the Bruhat ordering in
$W^\lambda$ is still generated by simple reflections,
\cite[Thm.14]{Caldero-Littelmann} also applies and yields that the
dual canonical basis of $V(n\lambda^*)$ coincides (up to nonzero
scalar multiples) with its standard monomial basis. As above, it
follows that any integral point of $Q_{\uw_0}(n\lambda)$ is a sum of
$n$ integral points of $Q_{\uw_0}(\lambda)$. This implies, in turn,
that the polytope $Q_{\uw_0}(\lambda)$ is integral.

It remains to show that $(G/P_\lambda)_0$ is a Fano. Note that 
$\cO(-K_{G/P_\lambda}) = L_{m\lambda}$ for some positive integer $m$
(as the weight $\lambda$ is fundamental), so that the polytope 
$Q_{\uw_0}(G/P_\lambda,\cO(-K))=mQ_{\uw_0}(\lambda)$ is integral. Now
Theorem \ref{QFano} completes the proof.
\end{proof}

%%%%%%%%%%%%%%%%%%%%%%%%%%%%%%%%%%%%%%%%
\section{Examples and counterexamples}
\label{sec:Examples and counterexamples}

%%%%%%%%%%%%%%%%%%%%
\subsection{$G$ of type $A_n$ and $\uw_0^{\rm std}$ (Gelfand-Tsetlin
  case)} 

The simplest reduced decomposition of the longest element in
$W=S_{n+1}$ is 
\begin{displaymath}
\uw_0^{\rm std}=(s_1)(s_2s_1)(s_3s_2s_1)\dots(s_ns_{n-1}\dots
s_1),
\end{displaymath}
where $s_i$ denotes the transposition exchanging $i$ with $i+1$.
It is convenient to use coordinates $x_{i,j}$ with $i,j\ge 1$,
$i+j\le n+1$ in place of $t_1, \dots, t_{n(n+1)/2}$. The string cone
is defined by
$$ 
x_{n,1}\ge0; \quad x_{n-1,2}\ge x_{n-1,1}\ge 0; \quad \dots \quad
x_{1,n} \ge \dots \ge x_{1,1} \ge 0,
$$
and the remaining $\lambda$-inequalities are
\begin{displaymath}
x_{i,j} \le \langle\lambda,\alpha_j^{\vee} \rangle - x_{i,j-1} +
\sum_{k=1}^{i-1}  (-x_{k,j-1} +2x_{k,j} - x_{k,j+1})
\end{displaymath}
for all $k=1,\dots, n$, $i=1,\dots, k$ in which we set
a variable $x$ to zero whenever its indices are out of bounds.

The corresponding polytopes $\Qwstd(\lambda)$ are known as
\emph{Gelfand-Tsetlin polytopes}, and we will use the notation
$\GT(\lambda)$ for them. A more familiar description of $\GT(\lambda)$
is in terms of \emph{Gelfand-Tsetlin patterns}
$$
g_{i,j} \ge g_{i+1,j} \ge g_{i,j+1} 
$$
in variables $g_{i,j}$, $i,j\ge 1$, $i+j\le n+1$ and with
variables $g_{0,i}$ corresponding to $\lambda$: 
$$
\lambda = \sum_{i=0}^{n+1} g_{0,i} \epsilon_i
= \sum_{i=0}^{n+1} \lambda_{i} \epsilon_i
, \quad
\text{so that} \quad 
\langle\lambda,\alpha_j^{\vee} \rangle = 
\langle\lambda,\epsilon_j - \epsilon_{j+1} \rangle = 
\lambda_j - \lambda_{j+1}.
$$
The inequalities can be represented by a graph $\Gamma$ in which
an arrow goes from a smaller variable to a larger:

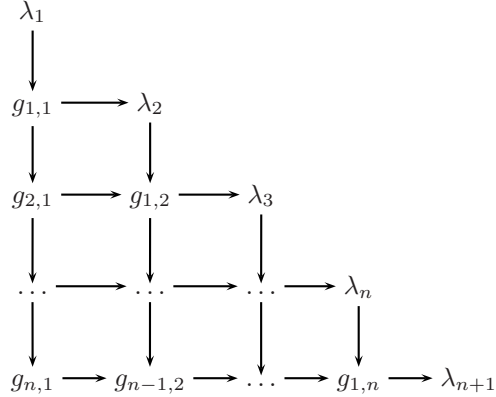
\begin{displaymath}
  \begin{psmatrix}[colsep=0.8cm,rowsep=0.8cm]
    \lambda_1 \\
    g_{1,1} & \lambda_2 \\
    g_{2,1} & g_{1,2} & \lambda_3 \\
    \dots   & \dots   & \dots & \lambda_n \\
    g_{n,1} & g_{n-1,2} & \dots & g_{1,n} & \lambda_{n+1} \\
    \everypsbox{\scriptstyle}
    \psset{arrows=->,nodesep=3pt}
    \ncline{1,1}{2,1}
    \ncline{2,1}{2,2}
    \ncline{2,1}{3,1}
    \ncline{2,2}{3,2}
    \ncline{3,1}{3,2}
    \ncline{3,1}{4,1}
    \ncline{3,2}{4,2}
    \ncline{3,2}{3,3}
    \ncline{3,3}{4,3}
    \ncline{4,1}{4,2}
    \ncline{4,1}{5,1}
    \ncline{4,2}{5,2}
    \ncline{4,2}{4,3}
    \ncline{4,3}{5,3}
    \ncline{4,3}{4,4}
    \ncline{4,4}{5,4}
    \ncline{5,1}{5,2}
    \ncline{5,2}{5,3}
    \ncline{5,3}{5,4}
    \ncline{5,4}{5,5}
  \end{psmatrix}
\end{displaymath}

The linear change of coordinates is given by
\begin{eqnarray*}
  x_{i,j} = \sum_{k=1}^j (g_{i-1,k} - g_{i,k}), \quad
  g_{i,j} = \lambda_j + \sum_{k=1}^i (x_{k,j-1} - x_{k,j} ).
\end{eqnarray*}

\begin{lemma}
  The fan $\Sigma_{\uwn^{\rm std}}$ is trivial. 
\end{lemma}
\begin{proof}
  Indeed, the additive property $\GT(\lambda+\mu) = \GT(\lambda) +
  \GT(\mu)$ is obvious from the Gelfand-Tsetlin patterns. (This was
  already noticed by Kaveh in \cite{Kaveh_SAGBI}.) 
\end{proof}

\begin{lemma}
  $\GT(\lambda)$ is integral $\iff$ 
  $\langle\lambda,  \alpha_{i}^{\vee} \rangle \in \bZ$ for all $i$.
\end{lemma}
\begin{proof}
  The vertices are given by setting some of the Gelfand-Tsetlin
  inequalities to equalities. It is then clear that $g_{i,j}$ are
  going to be equal to some $\lambda_j$'s and that all of
  them appear. So if the $\lambda_j$'s are integral
  then so are the $g_{i,j}$'s and vice versa.
\end{proof}

In fact, the following explicit description of the vertices is obvious:

\begin{lemma}
  Each vertex of the Gelfand-Tsetlin polytope $\GT(\lambda)$
  corresponds to setting each variable $g_{i,j}$ to some
  $\lambda_j$ in the way allowed by inequalities.  Equivalently,
  the boundaries between different values of $g_{i,j}$'s are paths
  $\Pi_i$ from the boundary of the graph $\Gamma$ along the diagonal
  to the lower left corner such that the union $\cup\Pi_i$ is a tree.
\end{lemma}

The collections $\cup \Pi_i$ were called \emph{meanders} in
\cite{Batyrev-et-al} in which the authors describe the toric
degeneration constructed by Gonciulea and Lakshmibai. Note that in the
case when the weight $\lambda$ is not regular, i.e. when some
$\lambda_j =0$ and the corresponding variety $G/P_{\lambda}$ is a
partial flag variety, the graph $\Gamma$ can be simplified by cutting
out square corners along the diagonal in which the values of
$g_{i,j}$'s are uniquely determined.

The extremal weight vertices of $\GT(\lambda)$ can be described with a
little more work. Since we do not need this description, we leave the
proof to the reader.

\begin{lemma}
  For every permutation $w\in W= S_{n+1}$, $w=(k_1 k_2 \dots k_{n+1})$ 
  the corresponding vertex is obtained by filling $k_1-1$ variables
  $g_{i,j}$ starting from the top in a minimal way with $\lambda_1$,
  then filling $k_2-1$ variables in a minimal way with
  $\lambda_2$ etc.
\end{lemma}

\begin{example}
  For the group $G=\SL_3$ of type $A_2$ and a regular weight
  $\lambda$, there are $6=3!$ extremal weight vertices and one
  additional vertex corresponding to setting
  $g_{1,1}=g_{1,2}=g_{2,1}=\lambda_2$.  The polytope
  $\GT(\lambda)$ has its projection to $\Lambda_{\bR}$ shown on
  Fig.~1.
\begin{figure}[htbp]
\label{fig:proj}
  \begin{center}
    \begin{pspicture}(-2.5,-2.1)(2.5,2.1)
      \definecolor{mygray}{gray}{0.8}
      \SpecialCoor
 %     \psgrid[gridcolor=Peach,subgridcolor=Apricot]
      \psdot(2;45)
      \uput[r](2;45){$\lambda$}
      \pspolygon[fillstyle=solid,fillcolor=mygray]%
      (2;45)(2;75)(2;165)(2;-165)(2;-75)(2;-45)
      \psline(2;75)(2;-75)
      \psline(2;165)(2;-45)
      \psline[linestyle=dashed](2;-165)(2;45)
      \psline[linecolor=gray](-2.5,0)(2.5,0)
      \psline[linecolor=gray](2.5;60)(-2.5;60)
      \psline[linecolor=gray](2.5;120)(-2.5;120)
    \end{pspicture}
  \end{center}
\caption{Projection of $\GT(\lambda)$}
\end{figure}
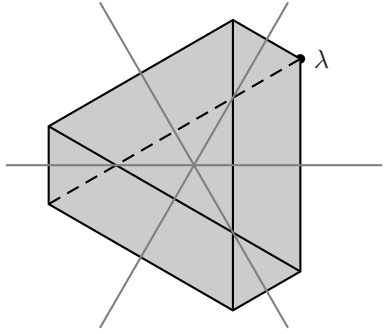

\end{example}

%%%%%%%%%%%%%%%%%%%%
\subsection{$G$ of type $A_n$ and other $\uw_0$}

For many of the explicit computations below, we used the freely
available program {\tt porta} for performing manipulations with
polytopes. 

The fans $\Sigma_{\uwn}$ need not be trivial in general. 

\begin{example}
  Let $G$ be of type $A_3$ and let the reduced decomposition be 
  $\uwn = s_1s_3s_2s_3s_1s_2$. The fan $\Sigma_{\uwn}$ consists of two
  maximal dimensional cones obtained by splitting $\Lambda^+_{\bR}$
  into two halves by the hyperplane 
  $\langle\lambda,\alpha_1^{\vee} \rangle =
  \langle\lambda,\alpha_3^{\vee} \rangle$. For a regular weight
  $\lambda$, the polytope $\Qlam$ has 38 or 44 vertices depending on
  which cone $\lambda$ lies in.
\end{example}

\begin{example}
  The polytopes $\Qw(\rho)$ are integral for $G= SL_n$, $n\le 5$.
  In type $A_3$ they have 12 or 13 facets and 38, 40 or 42
  vertices. In type $A_4$, the polytopes have from 20 to 27 facets and
  from 334 to 425 vertices.
\end{example}

\begin{conjecture}
  For $G$ of type $A_n$ and any reduced decomposition $\uw_0$, the
  polytope $Q_{\uw_0}(\lambda)$ is integral if and only if
  $\langle\lambda,\alpha_i^{\vee} \rangle \in\bZ$ for all $i$.
\end{conjecture}

%%%%%%%%%%%%%%%%%%%%
\subsection{Polytopes for fundamental weights}

\begin{example}
  A direct computation shows that for a group of type $A_n$, $n\le 4$,
  the polytopes $Q_{\uwn}(\omega_i)$ for a fixed fundamental weight
  $\omega_i$ are isomorphic, i.e. they do not depend on $\uwn$.
  However, they are different as polytopes marked by two distinguished
  vertices (the highest and lowest weight vertices). For example, in
  type $A_3$ the tangent cone at the origin of the polytope
  $\Qwstd(\omega_2)$ is non-simplicial. On the other hand, for 
  $\uwn = s_2 s_1 s_2 s_3 s_2 s_1$ the tangent cone of $\Qw(\omega_2)$
  at the origin is simplicial. 
\end{example}
%%%%%%%%%%%%%%%%%%%%

\subsection{Other classical types}

\begin{example}
  For $G$ of type $C_2$, Littelmann \cite[\S2]{Littelmann} computes
  the string cones for the reduced decompositions $\uw_0=s_1s_2s_1s_2$
  and $\uw'_0=s_2s_1s_2s_1$. They are respectively
  $$
  \{t_1\ge 0; \quad 2t_2 \ge t_3 \ge 2t_4 \ge 0\} \quad\text{and}
  \quad \{t_1\ge0; \quad t_2 \ge t_3 \ge t_4 \ge 0\}. $$
  In both cases, the fan $\Sigma$ is trivial. A simple
  computation shows that in the second case all the polytopes
  $Q(\lambda)$ are integral. However, in the first case,  
  the polytope $Q(\lambda)$ is integral if and only if
  $\langle\lambda,\alpha^{\vee}_1\rangle$ is even. In particular,
  $Q(2\rho) = Q(G/B, \cO(-K))$ is integral.  The polytopes
  $Q_{\uw_0}(2\rho)$ and $Q_{\uw'_0}(2\rho)$ both have 12 vertices.
\end{example}

There is a notion of Gelfand-Tsetlin polytopes in types $B$, $C$
and $D$. Their defining inequalities can be found in
\cite{BerensteinZelevinsky_JGP} and \cite[\S6,\S7]{Littelmann}. From
the definitions, it is obvious that in types $B_n$ and $C_n$ the fan
$\Sigma_{\uwn^{\rm std}}$ is trivial and the polytope
$\Qwstd(\lambda)$ is integral whenever the
$\langle\lambda,\alpha^{\vee}_i\rangle$ are all even. In particular,
the polytope $\Qwstd(2\rho)$ for the anticanonical limit of $G/B$ is
(integral and) reflexive, and the corresponding variety $(G/B)_0$ is 
a Fano.

In type $D_n$, the fan $\Sigma_{\uwn^{\rm std}}$ is non-trivial.

%%%%%%%%%%%%%%%%%%%%
\subsection{$G$ of exceptional type}

\begin{example}\label{ex:E6}
  For $G$ of type $E_6$ and some nice decomposition $\uw_0$, the polytope
  $Q_{\uw_0}(2\rho) = Q_{\uw_0}(G/B, \cO(-K))$ is not integral.
  
  Littelmann \cite[\S8]{Littelmann} computed (and we checked) the
  string cone for the reduced decomposition $\uw_0=\tau \cdot
  s_6s_2s_3 s_1s_4s_5 s_3s_4s_2 s_3s_1s_6 s_2s_3s_4 s_5$ with the
  enumeration 
  
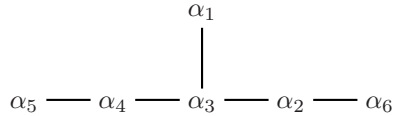
\begin{displaymath}
    \begin{psmatrix}[colsep=0.8cm,rowsep=0.8cm]
      && \alpha_1 \\
      \alpha_5 & \alpha_4 & \alpha_3 & \alpha_2 & \alpha_6 
      \everypsbox{\scriptstyle}
      \psset{arrows=-,nodesep=3pt}
      \ncline{1,3}{2,3}
      \ncline{2,1}{2,2}
      \ncline{2,2}{2,3}
      \ncline{2,3}{2,4}
      \ncline{2,4}{2,5}
    \end{psmatrix}
  \end{displaymath}
  of the simple roots, where $\tau$ is a reduced decomposition for the
  longest element in the Weyl group of the sub-root system of type $D_5$
  obtained by suppressing $\alpha_6$. The cone $C_{\uw_0}$ is the
  product of the string cone for $D_5$ and the cone

  \begin{eqnarray*}
  & t_1 \ge t_2 \ge t_3 \ge 
  \myfrac{t_4}{t_5} 
  \ge
  t_7 \ge 
  \myfrac{t_8}{t_9}
  \ge t_{10} \ge
  \myfrac{t_{11}}{t_{13}} 
  \ge 
  t_{14} \ge t_{15} \ge t_{16} \ge 0;  \\
  & t_5 \ge t_6 \ge t_8\ge 0; \quad
  t_9 \ge t_{12} \ge t_{13} \ge0.
  \end{eqnarray*}
  Hence, the projection of the 36-dimensional polytope
  $Q_{\uw_0}(\lambda)$ to the last 16 coordinates is the
  16-dimensional polytope $Q'$ defined by the inequalities above and
  the $\lambda$-inequalities for $k=36,35, \dots 21$. An explicit
  computation shows that for $\lambda=\rho$, i.e. when all
  $\langle\lambda,\alpha_j^{\vee} \rangle = 1$, the polytope $nQ'$ is
  integral if and only if $n$ is divisible by $6$. Hence, the
  polytopes $2Q'$ and $Q_{\uw_0}(2\rho)$ are not integral, and the
  limit toric variety of the anticanonically polarized variety $G/B$
  is not a Fano.
\end{example}

%%%%%%%%%%%%%%%%%%%%
\subsection{String polytopes and Duistermaat-Heckman measure}

An important question is what the string polytopes $\Qw(\lambda)$ have
in common for different reduced decompositions $\uwn$. Here, we note
two such properties: the Ehrhart function 
\begin{displaymath}
n\mapsto \#(\Qw(n\lambda)\cap \bZ^N) = \dim V(n\lambda^*)
= \dim V(n\lambda)
\end{displaymath}
for $n\ge0$, which is a polynomial function even if
$Q_{\uw_0}(\lambda)$ is not integral, and the Duistermaat-Heckman
measure.
  
Recall that for any symplectic (differentiable) manifold $(X,\omega)$
of dimension $2n$, equipped with a Hamiltonian action of a compact
torus $K$, the Duistermaat-Heckman measure $d\sigma_{DH}^K$ is the
pushforward of the Liouville measure $\omega^n/n!$ under the moment
map $\mu: X\to (\Lie K)^*$ to the dual of the Lie algebra of $K$. The 
image $\mu(X)$ is a convex polytope, and the measure
$d\sigma_{DH}^K$ is continuous and piecewise polynomial on this
polytope.
  
On the other hand, let $(X,L)$ be a (possibly singular) polarized
projective toric variety under the torus $\bT$. Denote by 
$\bK\subset \bT$ the maximal compact subtorus. Then $X$ admits a
moment map $\mu_{\bK}: X\to (\Lie \bK)^*$, which may be defined as
the composition of the finite morphism $X\to \bP H^0(X,L)^*$, with the
moment map for the $\bK$-action on $\bP H^0(X,L)^*$ (regarded as
a symplectic variety via a $\bK$-invariant K\"ahler structure). The
image of the moment map identifies with the moment polytope $P(X,L)$,
and the measure $d\sigma^{\bK}_{DH}$ is just the standard Euclidean
measure multiplied by $(2\pi)^{\dim X}$ on this polytope.

Given another torus $T$ with maximal compact subtorus $K$, and a
homomorphism $h:T\to\bT$, the $K$-action on $X$ via $h$ admits a
moment map $\mu_{K}: X\to (\Lie K)^*$ which is nothing but 
the map $\mu_{\bK}: X\to (\Lie \bK)^*$
followed by the map $\pi:(\Lie \bK)^* \to (\Lie K)^*$ dual of the 
map $\Lie K \to \Lie \bK$ induced by $h$. The corresponding
Duistermaat-Heckman measure $d\sigma^{K}_{DH}$ coincides, up to a
constant, with the pushforward measure of the linear map of
polytopes $\pi:\mu_{\bK}(X) \to \mu_{K}(X)$. It is piecewise
polynomial of degree $\dim X - \dim h(T)$.

In our case, the limit of the partial flag variety $G/P_\lambda$ has
two torus actions: by $\bT$ and by the maximal torus $T$ of $G$,
related by the homomorphism
\begin{displaymath}
h:T\to\bT, \quad t\mapsto (\alpha_{i_1}(t),\ldots,\alpha_{i_N}(t))
\end{displaymath}
with kernel being the center of $G$. The corresponding moment
polytopes are $\Qw(\lambda)$ and $\Conv(W\lambda^*)$.

\begin{lemma}
  The pushforward of the Euclidean measure under the projection
  $\pi_\lambda:\Qw(\lambda) \to \Conv (W\lambda^*)$ does not depend on
  $\uwn$.  
\end{lemma}
\begin{proof}
  Indeed, let $(\cX,\cO_{\cX}(1))\to \bA^1$ be a one-parameter
  degeneration with special fiber the toric variety corresponding to
  the polytope $\Qw(\lambda)$. The Duistermaat-Heckman measure
  depends continuously on $t\in \bC$ and is constant for $t\ne 0$. 
  Therefore, the measure in the limit coincides with the measure
  $d\sigma_{DH}$ for the $K$-action on $G/P_\lambda$ (regarded as
  the coadjoint orbit of $\lambda$ under a maximal compact subgroup of
  $G$).

  Alternatively, one may observe that the density of the pushforward
  of the Euclidean measure under $\pi_\lambda$ assigns to any 
  $\mu\in\Conv(W\lambda^*)$ the volume of its fiber
  $\pi_\lambda^{-1}(\mu)$. For rational $\mu$, this fiber is a
  rational convex polytope with Ehrhart function
  \begin{displaymath}
  n\mapsto \#(n\pi_{\lambda}^{-1}(\mu)\cap\bZ^N) = 
  \#(\pi_{n\lambda}^{-1}(n\mu)\cap\bZ^N) = 
  \dim V(n\lambda^*)_{n\mu},
  \end{displaymath}
  the multiplicity of the weight $n\mu$ in $V(n\lambda^*)$. Thus,
  the density function is independent of $\uw_0$.
\end{proof}

We note that this measure contains only partial information about
the polytope $\Qw(\lambda)$. For example, the measure on Fig.1 is
$S_3$-symmetric. However, the corresponding Gelfand-Tsetlin polytope
is not $S_3$-symmetric, and only one of the special interior points on
the picture is the image of a vertex of $\Qw(\lambda)$.

For any $\bQ$-polarized spherical variety $(X,L)$, the same proof
shows that the pushforward measure under the projection
$\pi:\Qw(X,L) \to \Conv(W\ P(X,L))$ does not depend on $\uwn$.

%%%%%%%%%%%%%%%%%%%%%%%%%%%%%%%%%%%%%%%%%%%%%%%%%%%%%%%%%%%
%%%%%%%%%%%%%%%%%%%%%%%%%%%%%%%%%%%%%%%%
\section{Connection with the ``universal'' compactified moduli space}
\label{Connection with the ``universal'' compactified moduli space}

The Minimal Model Program, which is fully established in dimension
$\le3$ and is still conjectural in dimensions $\ge4$, implies the
existence of a complete moduli space $\oM$ of stable pairs $(X,D)$,
generalizing the Deligne-Mumford moduli space $\overline{M}_g$ of
stable curves. These pairs consist of a projective variety $X$ and a
$\bQ$-Weil divisor $D$ that satisfy two basic conditions: 
\begin{enumerate}
\item (on singularities) $(X,D)$ should have semi-log canonical
  singularities, and
\item (numerical) the log canonical divisor $K_X + D$, properly defined
  if $X$ is non-normal, should be ample. 
\end{enumerate}
In particular, every one-parameter family $(X_t, D_t)$ of stable pairs
should have, possibly after a finite base change, a unique stable
limit $(X_0,D_0)$. 
See \cite{KollarShepherdBarron,Alexeev_LCanModuli,AlexeevMgn} for more
of this story. 

In the situations where the varieties have an additional structure,
such as a group action, the existence of such moduli space can be
established in higher dimensions without using the methods of the
Minimal Model Program. Two examples are pairs with toric or
semiabelian group action \cite{Alexeev_CMAV} and stable reductive
varieties \cite{AlexeevBrion_Affine,AlexeevBrion_Projective}.

Below, we show that in the cases when the polytope $Q_{\uwn}(X,L)$ is
integral the toric limits we have considered in this paper can be
interpreted as boundary points on the compactified moduli space $\oM$
of stable pairs, it it exists. First, recall the basic

\begin{definition}
  Let $X$ be a normal variety and $D=\sum a_j D_j$ be a $\bQ$-Weil
  divisor such that $0\le a_j \le 1$ for all $j$. The pair 
  $(X,D)$ is canonical if some positive multiple of $K_X + D$ is an
  integral Cartier divisor, and if for every resolution of
  singularities $f:Y\to X$ with exceptional divisors $F_i$, 
  we have
  \begin{displaymath}
  K_Y  + f\inv_{\rm strict} D = 
  f^*(K_X + D) + \sum b_i F_i
  \end{displaymath}
  where $b_i\ge -1$ for all $i$. This formula makes sense since
  Cartier divisors can be pulled back; the coefficients $b_i$ are
  called the discrepancies.
\end{definition}

The notion of semi-log canonical pair is a generalization of log
canonical to the case when $X$ is non-normal. We will not need this
generalization since all our varieties are normal.

Let $X=G/B$ be a flag variety, $E$ be the union of all the Schubert
varieties considered as a reduced divisor, and $E^-=w_0 E$ be the
opposite divisor. For a regular dominant weight $\lambda$, define the
divisor $D=D_{\lambda}$ on $X$ as follows:
\begin{displaymath}
  D= (s), \quad \text{ where} \quad s = 
  \sum  b_{\lambda,\varphi} \in H^0(G/B,L_\lambda) 
\end{displaymath}
with the latter sum going over all elements of the dual canonical
basis of $V(\lambda^*)$.

\begin{theorem}
\begin{enumerate}
\item 
The pair $(X,E + E^-)$ is log canonical and $K_X + E + E^-$ is
rationally equivalent to $0$.
\item 
For $0<\varepsilon \ll1$, the pair $(X,E + E^- + \varepsilon D)$
    is stable, i.e. $(X,E + E^- + \varepsilon D)$ is log canonical and
    $K_X+ E + E^- + \varepsilon D$ is ample. 
\end{enumerate}
\end{theorem}
\begin{proof}
  (i) is proved in \cite[Thm.5.10]{AlexeevBrion_Projective}. Since $D$
  is ample then the divisor $K_X+E+E^-+\varepsilon D$ obviously is
  ample. To prove the statement about log canonical singularities, we
  have to prove that $D$ does not contain any center of log canonical
  singularities of $E+E^-$, i.e. any image of an exceptional divisor
  $F_i$ on a resolution $f:Y\to X$ that has discrepancy $-1$. But
  such a center must be $T$-invariant, and $s$ does not vanish
  at any $T$-fixed point of $X$.
\end{proof}

Now, pick any family of integers $(n_{\varphi})$. Then the section
\begin{displaymath}
s_t = \sum t^{n_{\varphi}} b_{\lambda,\varphi} 
\in H^0(G/B, L_{\lambda})
\end{displaymath}
defines a one-parameter family of stable pairs $(X,E+E^-+\varepsilon
D_t)$. The Minimal Model Program in dimension $\dim(X)+1$ implies that
this family has a unique limit as a stable pair
$(X_0,E_0+E_0^-+\varepsilon D_0)$. It is obtained as the special fiber
of the log canonical model $\cX_{\rm can}$ of the semistable model of
the degenerating family $\cX$.

On the other hand, the filtration $(R_{\le m})$ of the algebra 
$R=R(G/B,L_{\lambda})$, constructed in
Proposition~\ref{deformation}, can be interpreted as the data for a
one parameter degeneration, in the usual way, by setting
$$ n_{\varphi} = \min\{ m \,|\, b_{\lambda,\varphi} \in R_{\le m} \}. $$

Now assume that the polytope $\Qw(\lambda)$ is integral.  On the limit
toric variety $X_0$, the limit of $E+E^-$ is the boundary divisor
$\partial X_0$, i.e., the complement of the open $\bT$-orbit.

Indeed, on $X = G/B$ the divisor $E+E^-$ is the divisor of zeros of
the product of the highest and the lowest weight vectors in
$H^0(G/B,L_\lambda)$ (since $\lambda$ is regular).  Thus, $E_0 +
E_0^-$ is the divisor of zeroes of the product of the limits of these
vectors on $X_0$. These limits are the sections associated with the
highest (resp. lowest) weight vertices of the string polytope. The
product of these sections vanishes everywhere on the boundary of
$X_0$, since any codimension-1 face of the string polytope contains
one of these vertices, as is clear from the defining inequalities
in Theorem~\ref{thm:cone_eqns}.

This proves that $\supp(E_0+E_0^-) = \supp \partial X_0$. As on any
toric variety, one has $K_{X_0}+\partial X_0 \sim 0$. Applying the
argument of Proposition \ref{singularities}, we see that $K_{\cX}+ \cE
+ \cE^- \sim 0$ on the family $\cX$. Therefore, $E_0+E_0^- \sim
\partial X_0 $.

This implies that the pair $(X_0, E_0+E_0^-)$ has log canonical
singularities and that its centers of log canonical singularities are
precisely the closures of $\bT$-orbits. The limit of $s_t$ is the
section $s_0=\sum b_{\lambda,\varphi}$, with $b_{\lambda,\varphi}$
going over the basis of $\bT$-eigenfunctions of $H^0(X_0,L_0)$. Again,
the divisor of this section does not contain any $\bT$-fixed points.
Hence, the pair $(X_0,E_0+ E^-_0 + \varepsilon D_0)$ is stable, as
predicted.

This argument can be repeated for an arbitrary polarized projective
spherical variety $(X,L)$, provided the polytope $Q_{\uw_0}(X,L)$ is
integral, with $s$ being replaced by the sum of basis vectors
$b_{\lambda,\varphi}$ with $\lambda$ going over the weights appearing
in $H^0(X,L)$. This concludes our speculation about the connection
with the ``universal'' moduli space of stable pairs. 

Finally, we want to remark on the construction of Grossberg and
Karshon \cite{GrossbergKarshon} of toric limits of Bott-Samelson
varieties. The latter appear naturally as resolutions of singularities
of the pair $(G/B, E)$, where $E$, as above, is the union of Schubert
varieties, considered as a reduced divisor. Every reduced
decomposition $\uwn$ provides a Bott-Samelson variety $Y=Y_{\uwn}$
together with a birational morphism $f=f_{\uwn}:Y \to G/B$.

The Bott-Samelson varieties degenerate to toric limits for the torus
$\bT$, that are are nonsingular toric varieties called Bott towers. In
\cite{GrossbergKarshon} they are described by a combinatorial
structure called \emph{twisted cube} which is not a polytope but a
kind of ``virtual'' polytope. Translating this back to geometry,
this means that the limit on $Y_0$ of the globally generated
line bundle $f^*(L_\lambda)$ is a sheaf which is no longer
globally generated. Hence, the relationship between $G/B$ and
$Y$ gets broken in the limit, and there is only a rational map
between the special fibers of degenerations.

%\bibliography{tordeg}
\providecommand{\bysame}{\leavevmode\hbox to3em{\hrulefill}\thinspace}
\providecommand{\MR}{\relax\ifhmode\unskip\space\fi MR }
\providecommand{\MRhref}[2]{%
  \href{http://www.ams.org/mathscinet-getitem?mr=#1}{#2}
}
\providecommand{\href}[2]{#2}

\end{document}